\documentclass[12pt]{amsart}
 \usepackage{amsfonts,amssymb,eucal}
\usepackage{amsthm} 
 \usepackage{amsmath} 
\usepackage{amscd}
 \usepackage{latexsym} 
\input{cyracc.def} 

\numberwithin{equation}{section}

\usepackage{color}
\newcommand{\dd}{\operatorname{d}}

\newcommand{\e}{{\mbox{\rm e}}}

\newcommand{\mb}[1]{{\mbox{\boldmath{$#1$}}}}
\newcommand{\mc}[1]{{\mathcal{#1}}}
\newcommand{\mr}[1]{{\mathrm{#1}}}
\newcommand{\got}[1]{{\mathfrak{#1}}}
\newcommand{\db}[1]{{\mathbb{#1}}}

\newcommand{\pa}{\partial}

\newcommand{\R}{\ensuremath{\mathbb{R}}}
\newcommand{\C}{\ensuremath{\mathbb{C}}}
\newcommand{\N}{\ensuremath{\mathbb{N}}}

\newcommand{\Hi}{\ensuremath{\got{H}}}
\def\ii{\operatorname{i}}
\newcommand{\tr}{\ensuremath{{\mbox{\rm{Tr}}}}}%
\newcommand{\zn}{\ensuremath{\mathbb{O}_n}}
\newcommand{\h}{\ensuremath{\got{h}}}

\newtheorem{Theorem}{Theorem}
\newcommand{\un}{\ensuremath{\mathbb{I}_n}}
\newtheorem{Proposition}{Proposition}

\theoremstyle{definition}

\begin{document}
\title{Geodesics associated to the balanced metric on the Siegel-Jacobi ball }
\author{Stefan  Berceanu}
\address[Stefan  Berceanu]{National Institute for Physics and Nuclear Engineering\\
         Department of Theoretical Physics\\
         PO BOX MG-6, Bucharest-Magurele, Romania}
\email{Berceanu@theory.nipne.ro}

\begin{abstract}
We determine the Christoffel's symbols for the Siegel-Jacobi ball
endowed with the balanced metric.  We study the equations of geodesics on
the Siegel-Jacobi ball. We calculate the covariant derivative of
one-forms in the variables in which is expressed the balanced metric
on the Siegel-Jacobi ball.
\end{abstract}
\maketitle
\noindent
\tableofcontents

\section{Introduction}\label{intro}
The Jacobi group of index $n$, $G^J_n$, defined as the semidirect product of the symplectic
group with Heisenberg group of appropriate dimension, is an interesting  object in several 
branches of Mathematics and Physics, see references  in
\cite{jac1,sbj,BERC08B,gem,SB15}. $G^J_n$ acts transitively   on the homogeneous
manifold associated with it, the Siegel-Jacobi ball $\mc{D}^J_n$, a
partially bounded domain, 
whose points are in $\C^n\times \mc{D}_n$, where $\mc{D}_n$ denotes
the Siegel ball. In fact, there is a real  and a complex version of the
Jacobi group, with  factor  group $\text{Sp}(n,\R)$,
respectively $\text{Sp}(n,\R)_{\C}$, while the associated manifold in
the real case is the Siegel-Jacobi upper half plane $\mc{X}^J_n$, whose
points are in $\C^n\times\mc{X}_n$, where $\mc{X}_n$ is the
Siegel upper half plane. The case of the Jacobi group of index $n=1$
was considered in \cite{ez,bs}.  There are generalizations of the
Jacobi group and their associated homogeneous spaces, where instead
$\C^n$  are considered   spaces $\C^{mn}$, $m,n\in\N$. There is an
isomorphism between the complex and real version of the Jacobi
group and the homogenous spaces associated to them \cite{Y10,nou}. The present paper is mainly devoted to the complex version of
the Jacobi group. 

In \cite{sbj,nou} we have determined the $G^J_n$-invariant metric on
$\mc{D}^J_n$ and we have shown that it corresponds to the metric
determined in the papers of J.-H. Yang \cite{Y10}, generalizing the results of
E. K\"ahler \cite{cal}  and R. Berndt \cite{ez}  in the case $n=1$. We
mention that the starting point of our approach  in \cite{sbj,nou}
was the construction of coherent states (CS) \cite{perG} attached to
the Jacobi group $G^J_n$, with support in $\mc{D}^J_n$. The case $n=1$
was studied in \cite{jac1}. What we have find starting from the
construction of CS attached to $G^J_n$ was the K\"ahler potential on
$\mc{D}^J_n$. In \cite{SB15} we have underlined that the metric
determined in \cite{sbj,nou} is the  balanced metric \cite{don,arr} and
we have determined the metric matrix $h(z,W)$ associated to the K\"ahler
two-form $\omega _{\mc{D}^J_n} (z,W)$, where $z\in\C^n$ and 
$W\in M(n,\C)$ is a symmetric matrix describing points in
$\mc{D}_n$. We have introduced  \cite{nou} a  useful variable $\eta\in\C^n$, such
that if we make a change of variables $ (z,W)=FC(\eta,W)$, $z=\eta-W\bar{\eta}$, the
metric on the Siegel-Jacobi ball expressed in $(\eta,W)$ separates as sum
of the metric on $\C^n$ and $\mc{D}_n$, i.e. the $FC$-transform is a
K\"ahler homogeneous diffeomorphism and realizes the fundamental conjecture for
the homogeneous K\"ahler manifold $\mc{D}^J_n$  \cite{GV,DN}.  In
\cite{SB15} we have also determined the inverse matrix $h^{-1}$ in the
same variables. 

In the present paper  we use the results of \cite{SB15}  and calculate
the Christoffel symbols on $\mc{D}^J_n$. This allows us to 
write down explicitly the equations of geodesics on $\mc{D}^J_n$. We
are not able to integrate these equations, but we make some comments
  related to the geometric and physical  meaning of
the variables used. In order to clarify the physical meaning of the
$FC$-transform, we use some notions concerning the CS, even that the
present work is mainly  a simple matrix  calculation,  where the CS does
not appear  explicitely. As a byproduct
of the calculation of the $\Gamma$-s on $\mc{D}^J_n$, we also
calculate the covariant derivatives  $D(\dd z)$ and $ D(\dd W)$
(called in \cite{YY1} holomorphic $\Gamma$-module connections
associated to the Christoffel symbols) - this kind of 
calculation was done also in \cite{YY2} in order to determine
derivations of the Jacobi forms on
$\mc{X}_{mn}=\mc{X}_n\times\C^{mn}$. The reason to present the
calculation of $D(\dd z)$ and $ D(\dd W)$ is to underline again the
utility of the variables introduced in \cite{sbj,nou}.

The paper is divided in two parts. In \S \ref{PR} are recalled some
notations and concepts introduced in pervious works  and used in the
present paper. In \S
\ref{JCG} we present what we mean by Jacobi group in this paper, see
details in \cite{sbj,nou}. \S \ref{BMCS} recalls some simple facts
about the CS defined 
by Perelomov \cite{perG}, in order to understand how the balanced
metric was obtained in \cite{SB15}. The formulae used for calculation
of geodesics are
recalled in \S \ref{GD}.  The geometry of $\mc{D}^J_1$
\cite{jac1,csg},  summarized  in \S
\ref{BSG}, is a guide for the calculation to be done on
$\mc{D}^J_n$, especially Proposition \ref{GER}, which gives equations
of geodesics on the Siegel-Jacobi disk.  \S \ref{BMSB} reproduces the results from \cite{SB15}
concerning the  matrix of the balanced metric on $\mc{D}^J_n$ and is its inverse, the starting point of
the present investigation - see Theorem \ref{mainTH}. In Proposition
\ref{TOTU} we extract from \cite{SB15} some facts on  the geometry of the
Jacobi group and its action on the Siegel-Jacobi ball. The new results
of the present paper are contained in \S \ref{P2}. The
equations of geodesics in the variable $z\in\C^n$, $W\in\mc{D}_n$ are
given  in \S \ref{S1}, which also contains the calculation of $\Gamma$-s. Equations of geodesics are collected in
Proposition \ref{PR2} of \S \ref{S3}, followed by a discussion about
the significance of the result in the context of the $FC$-transform on
$\mc{D}^J_n$,
which is shown to be not a geodesic mapping. It was proved
\cite{ber97} that for symmetric manifolds the  $FC$-transform gives
geodesics, but the studied Siegel-Jacobi ball is not a symmetric space.  Section \ref{CDDD}
summarizes our calculation of the covariant derivative of $\dd z$ and
$\dd W$.  The  appendix in \S \ref{APP} is dedicated to the equations of geodesics
on $\mc{D}_n$ deduced from equations of geodesics on $\mc{X}_n$. The
new results of these paper are contained in Propositions \ref{CRR},
\ref{PR2}, \ref{part} and \ref{last}.

\textbf{Notation.} 
We denote by $\mathbb{R}$, $\mathbb{C}$, $\mathbb{Z}$,
and $\mathbb{N}$ the field of real numbers, the field of complex numbers,
the ring of integers, and the set of non-negative integers, respectively.  We denote the imaginary unit
$\sqrt{-1}$ by $\ii$, and the Real and Imaginary part of a complex
number by $\Re$ and respectively $\Im$, i.e. we have for $z\in\C$,
$z=\Re z+\ii \Im z$, and $\bar{z}=  \Re z-\ii \Im z$.
$M(m\times n, \mathbb{F})\approxeq\mathbb{F}^{mn}$ denotes the set of all $m\times
n $ matrices with entries in the field $\mathbb{F}$. $M(n\times 1,\mathbb{F})$ is
identified with $\mathbb{F}^n$. Set $M(n,\mathbb{F})=M(n\times n,\mathbb{F})$.
For any $A\in M_{n}(\mathbb{F})$, $A^{t}$ denotes the transpose matrix of
$A$. 
 For $A\in M_{n}(\mathbb{C})$, $\bar{A}$ denotes the conjugate matrix
of $A$ and $A^{*}=\bar{A}^{t}$. For $A\in M_n(\mathbb{C})$, the
inequality $A>0$ means that $A$ is positive definite. The identity matrix of
degree $n$ is denoted by $\un$ and  $\zn$ denotes the $M(n,\mathbb{F})$-matrix
with all entries $0$. 
We consider a complex separable Hilbert space $\got{H}$ endowed with a   scalar product
which  is antilinear in the first argument,
$(\lambda x,y)=\bar{\lambda}(x,y)$, $x,y\in\got{H}$,
$\lambda\in\C\setminus 0$. If $A$ is a linear operator, we denote by
$A^{\dagger}$ its adjoint. If $\pi$ is a  representation of a Lie group
$G$ on the Hilbert space \Hi~ and $\got{g}$ is the Lie algebra of $G$,
we denote $\mb{X}:=\dd \pi(X)$  for $X\in\got{g}$. 
 We use Einstein convention  that repeated indices are
implicitly summed over.
\section{Preliminaries}\label{PR}
\subsection{The Jacobi group}\label{JCG}

 The (complex) Jacobi group of index $n$ is defined
 as the semidirect product 
$G^J_n=\mr{H}_n\rtimes\text{Sp}(n,\R)_{\C}$,   where $\mr{H}_n$ denotes the 
$(2n+1)$-dimensional Heisenberg group \cite{Y02,sbj,nou}, endowed with  
the composition law 
\begin{equation}\label{compositieX}
(g_1,\alpha_1,t_1)(g_2,\alpha_2, t_2)= \big(g_1  g_2,
g_2^{-1}\times \alpha_1+\alpha_2, t_1+ t_2 +\Im
(g^{-1}_2\times\alpha_1\bar{\alpha}_2)\big).
\end{equation}
$\alpha_i \in\C^n$, $t_i\in\R$ and  $g_i\in\text{Sp}(n,\R)_{\C}$, 
$i=1,2$  have the form
\eqref{dgM}
\begin{equation}\label{dgM}
g = \left( \begin{array}{cc}p & q\\ \bar{q} &
\bar{p}\end{array}\right),~p,q\in M(n,\C), 
\end{equation}
  and $g\times\alpha = p\,\alpha +
q\,\bar{\alpha} $, and $g^{-1}\times\alpha ={p}^*\alpha
-q^t\bar{\alpha}$.
The matrices $p,q\in M(n,\C)$ have the properties
\begin{subequations}\label{simplectic}
\begin{eqnarray}
& pp^*- qq^* = \un,\quad  pq^t=qp^t; \label{simp1}\\
& p^*p-q^t\bar{q} = \un ,\quad p^t\bar{q}=q^*p .\label{simp2}\end{eqnarray}
\end{subequations}

To the Jacobi
group $G^J_n $ it is associated  a homogeneous manifold - the
Siegel-Jacobi  ball  $\mc{D}^J_n$ \cite{sbj} - whose points are in
$\C^n\times\mc{D}_n$, i.e. a partially-bounded space.   $\mc{D}_n$ denotes the
Siegel (open)  ball  of index $n$. The non-compact hermitian
symmetric space $ \operatorname{Sp}(n, \R
)_{\C}/\operatorname{U}(n)$ admits a matrix realization  as a  homogeneous bounded
domain:
\begin{equation}\label{dn}
\mc{D}_n:=\{W\in  M (n, \C ): W=W^t, N>0, N:=\un-W\bar{W} \}.
\end{equation}
For  $g\in
\operatorname{Sp}(n,\R )_{\C}$   of the form \eqref{dgM}, \eqref{simplectic}
and $\alpha, z\in\C^n$,  the transitive action $(g,\alpha)\times(W,z)=(W_1,z_1)$ of the Jacobi group $G^J_n$ on the
Siegel-Jacobi ball $\mc{D}^J_n$  is given by the formulae  \cite{sbj}
\begin{subequations}\label{TOIU}
\begin{align}
\label{x44} W_1 & =(pW+q)(\bar{q}W+\bar{p})^{-1}=(Wq^*+p^*)^{-1}(q^t+Wp^t),\\
\label{xxxX}
z_1 & = (Wq^*+ p^*)^{-1}(z+ \alpha -W\bar{\alpha}).
\end{align}
\end{subequations}

\subsection{Balanced metric and coherent states}\label{BMCS}
We consider a K\"ahler manifold $M$ endowed with a K\"ahler-two form,
locally written as 
\begin{equation}\label{kall}
\omega_M(z)=\ii\sum_{\alpha,\beta=1}^n h_{\alpha\bar{\beta}} (z) \dd z_{\alpha}\wedge
\dd\bar{z}_{\beta}, ~h_{\alpha\bar{\beta}}= \bar{h}_{\beta\bar{\alpha}}= h_{\bar{\beta}\alpha},
\end{equation}

We consider a   $G$-invariant  K\"ahler two-form $\omega_M$ \eqref{kall} on the $2n$-dimensional
  homogeneous manifold  $M=G/H$  derived from the  K\"ahler potential $f(z,\bar{z})$ \cite{chern}
\begin{equation}\label{Ter} 
h_{\alpha\bar{\beta}}= \frac{\pa^2 f}{\pa {z}_{\alpha}\pa
  \bar{z}_{{\beta}}} .
\end{equation}
 The condition of the metric to be a K\"ahlerian one is (cf. (6) p. 156 in \cite{koba2})
\begin{equation}\label{condH}
\frac{\pa h_{\alpha\bar{\beta}}}{\pa z_{\gamma}} = \frac{\pa
  h_{\gamma\bar{\beta}}}{\pa z_{\alpha}},~\alpha,\beta,\gamma
=1,\dots, n .
\end{equation}
As was underlined in \cite{berr} for $\mc{D}^J_1$ and proved in
\cite{SB15} for $\mc{D}^J_n$ , the homogeneous  hermitian  metric determined in \cite{jac1,sbj,nou}
is in fact a balanced metric, because it corresponds  to the K\"ahler
potential calculated as the scalar  product of two CS-vectors 
\begin{equation}\label{FK}
f(z,\bar{z})=\ln K_M(z,\bar{z}), ~~K_M(z,\bar{z})=(e_{\bar{z}},e_{\bar{z}}).
\end{equation}
We recall that in the approach
of Perelomov \cite{perG} to 
 CS  it is supposed that there exists 
 a continuous, unitary, irreducible 
representation $\pi$
 of a   Lie group $G$
 on the   separable  complex  Hilbert space \Hi . 
  The coherent vector
 mapping $\varphi$ is defined locally
 (cf. \cite{last,sb6}) $
 \varphi : M\rightarrow \bar{\Hi}, ~
 \varphi(z)=e_{\bar{z}}$,
where $ \bar{\Hi}$ denotes the Hilbert space conjugate to $\Hi$.

We
can introduce the normalized (unnormalized) vectors  $\underline{e}_x$
(respectively, $e_z$) defined on $G/H$
\begin{equation}\label{bigch}
\underline{e}_x=\exp(\sum_{\phi\in\Delta_+}x_{\phi}\mb{X}_{\phi}^+-\bar{x}_{\phi}\mb{X}_{\phi}^-)
e_0, ~
e_z=\exp(\sum_{\phi\in\Delta_+}z_{\phi}\mb{X}_{\phi}^+)e_0, 
\end{equation}
where
$e_0$ is the extremal weight vector of the representation $\pi$, $\Delta_+$ are the positive roots of the Lie algebra $\got{g}$ of $G$,
and $X_\phi,\phi\in\Delta$,  are the  generators. 
$\mb{X}^+_{\phi}$ ($\mb{X}^-_{\phi}$)
corresponds to the positive (respectively, negative) generators. See
details in  \cite{perG,sb6}.

Let us denote by $FC$ the change of variables
$x\rightarrow z$ in formula \eqref{bigch} such that
\begin{equation}\label{etild}\underline{e}_{x}=\tilde{e}_z, ~    \tilde{e}_z  :=
  (e_z,e_z)^{-\frac{1}{2}}e_z, ~z=FC(x). 
\end{equation}
The reason for calling the transform \eqref{etild}  $FC$  {\it
  (fundamental conjecture)} is explained later, see Proposition \ref{expl}. We
also recall that we have proved that {\it  for symmetric spaces
the dependence} $z(t)=FC(tX)$ from \eqref{etild} {\it gives geodesics in} $M$
{\it with the property that} $z(0)=p$ {\it and} $\dot{z}(0)=X$ \cite{ber97}.

Using Perelomov's CS vectors
$e_{\bar{z}}\in\bar{\Hi}$,  $z \in M $  
we have considered \cite{jac1,sbj,nou,SB15} Berezin's approach to quantization on K\"ahler
manifolds, with the supercomplete set of vectors
verifying the Parceval   overcompletness identity \cite{ber73,ber74,ber75,berezin}
\begin{equation}\label{PAR}
(\psi_1,\psi_2)_{\mc{F}_{K}}=\int_M (\psi_1,e_{\bar{z}}) (e_{\bar{z}},\psi_2) \dd
\nu_M(z,\bar{z}), \quad \psi_1,\psi_2\in\Hi ,
\end{equation}
\begin{equation}\label{DELNU}
\dd{\nu}_{M}(z,\bar{z})=\frac{\Omega_M(z,\bar{z})}{(e_{\bar{z}},e_{\bar{z}})};
\quad \Omega_M:=\frac {1}{n!}\;
\underbrace{\omega\wedge\ldots\wedge\omega}_{\text{$n$ times}}.
\end{equation}

Let us denote by   $\Hi_f$  the
weighted Hilbert space of square integrable holomorphic functions on
$M$, with weight $\e^{-f}$  \cite{eng}
\begin{equation}\label{HIF}
\Hi_f=\left\{\phi\in\text{hol}(M) | \int_M\e^{-f}|\phi|^2 \Omega_M
  <\infty \right\}.
\end{equation}

In order to identify the Hilbert space $\Hi_f$ defined by \eqref{HIF}
with the Hilbert ${\mc{F}_{K}}$ space with scalar product \eqref{PAR},
it was introduced 
  the  so called $\epsilon$-function  \cite{raw,Cah, cah}
\begin{equation}\label{EPSF}
\epsilon(z) = \e^{-f(z)}K_M(z,\bar{z}).
\end{equation}
If the K\"ahler metric on the complex manifold  $M$ 
 is such that
$\epsilon(z)$ is a positive  constant, then the metric is called {\it balanced}.
Donaldson   \cite{don}   used this denomination  for compact manifolds, then
it was used in \cite{arr} for noncompact manifolds and later   \cite{alo} in
the context of Berezin quantization on homogeneous bounded
domains.

The {\it balanced hermitian  metric} of $M$ in
local coordinates  is 
\begin{equation}\label{herm}
\dd s^2_M(z,\bar{z}) =\sum_{\alpha,\beta=1}^n \frac{\pa^2}{\pa
  z_{\alpha} \pa\bar{z}_{\beta}} \ln (K_M(z,\bar{z})) 
\dd z_{\alpha}\otimes \dd\bar{z}_{\beta} ,
~K_M(z,\bar{z})=(e_{\bar{z}},e_{\bar{z}}). 
\end{equation}

\subsection{Geodesics on K\"ahler manifolds}\label{GD}

In terms of a local coordinate system $x^1,\dots$, $x^n$ the equations of geodesics on a manifold $M$ with linear connection
with components of the linear connections $\Gamma$ are  (see
e.g. Proposition 7.8 p. 145 in \cite{koba1}) 
\begin{equation} \label{GEOO}
\frac{\dd ^2 x^i}{\dd t^2}  +\sum_{j,k}\Gamma^i_{jk}\frac{\dd x^j}{\dd t}
\frac{\dd x^k}{\dd t}  =  0, ~i=1,\dots,n .
\end{equation}
The components $\Gamma^i_{jk}$ of a Riemannian (Levi-Civitta) connection 
are obtained by solving  the linear system  (see e.g. \cite{koba1}, p 160)
\begin{equation}\label{geot}h_{lk}\Gamma^l_{ji}=\frac{1}{2}\left(\frac{\pa
      h_{ki}}{\pa x_j}+\frac{\pa h_{jk}}{\pa x_i}-\frac{\pa h
      _{ji}}{\pa x_k}\right).\end{equation}
The connection matrix (form) $\theta=\pa \ln h$ (see \S 6 in \cite{chern} or
\S 3.2 in \cite{ball}) has the matrix elements:
\begin{equation}\label{CV1}
\theta^j_i=\Gamma^j_{ik}\dd x^k.\end{equation}
The covariant derivative of a contravariant vector (one-form) is given
by 
\begin{equation}\label{CV2}
 D u_i=-\theta^i_j u_j.\end{equation}
If we take into account the hermiticity condition  in \eqref{Ter}  on
the metric   and the  K\"ahlerian
restrictions \eqref{condH}, 
the non-zero Christoffel's  symbols  $\Gamma$  of the Chern connection (cf. e.g. \S 3.2 in \cite{ball};
in the  case of K\"ahler manifolds, also Levi-Civita connection,  cf. e.g. Theorem 4.17
in \cite{ball}) on K\"ahler  manifolds which appear in \eqref{geot} are determined by the equations
(see also  e.g. (12) at p. 156  in \cite{koba2})  
\begin{equation}\label{CRISTU}  h_{\alpha\bar{\epsilon}}\Gamma^{\alpha}_{\beta\gamma}=
\frac{\pa  h_{\bar{\epsilon}\beta}}{\pa z_{\gamma}} = \frac{\pa
  h_{\beta\bar{\epsilon}}}{\pa z_{\gamma}}. 
\end{equation}
We introduce the  convention
\begin{equation}\label{sum1}
h^{\alpha\bar{\beta}}:=(h_{\alpha\bar{\beta}})^{-1},~~\text{i.e.~~}
h_{\alpha\bar{\epsilon}}h^{\epsilon \bar{\beta}}=\delta_{\alpha\beta}.
\end{equation}
From \eqref{CRISTU}, we have:
$$h_{\alpha\bar{\epsilon}}h^{\epsilon\bar{\eta}}\Gamma^{\alpha}_{\beta\gamma}=
  \delta^{\eta}_{\alpha}\Gamma^{\alpha}_{\beta\gamma}=h^{\epsilon\bar{\eta}}\frac{\pa
    h_{\gamma\bar{\epsilon}}}{\pa z_{\beta}},$$
and we find for the Christoffel's symbols for the K\"ahler manifold, the  formula 
\begin{equation}\label{stefaN}
\Gamma^{\gamma}_{\alpha\beta}=\bar{h}^{\gamma\bar{\epsilon}}\frac{\pa
  h_{\beta\bar{\epsilon}}}{\pa z_{\alpha}}=
h^{\epsilon\bar{\gamma}}\frac{\pa h_{\beta\bar{\epsilon}}}{\pa  z_{\alpha}} .
\end{equation}

 \subsection{Geodesics on the Siegel-Jacobi disk}\label{BSG}

We recall some formulae describing the the Jacobi group $G^J_1$ and
the geometry of the 
Siegel-Jacobi  disk $\mc{D}^J_1$ \cite{jac1,csg,berr}.

The transitive action of the group $ G^J_1=\mr{H}_1\rtimes\text{SU}(1,1)\ni
(g,\alpha)\times(z,w)\rightarrow(z_1,w_1)\in \mc{D}^J_1$ on the Siegel disk
is given by the formulae:
\begin{equation}\label{dg}  w_1 
=\frac{a \, w+ b}{\delta}, ~ \delta=\bar{b}w+\bar{a},
\text{SU}(1,1) \ni g= \left( \begin{array}{cc}a & b\\ \bar{b} &
\bar{a}\end{array}\right),~\text{where} ~|a|^2-|b|^2=1,
\end{equation} 
\begin{equation}\label{xxx}
z_1=\frac{\gamma}{\delta}, ~\gamma = z
+\alpha-\bar{\alpha}w. 
\end{equation}
The balanced  K\"ahler two-form $\omega_{k\mu}$ on $\mc{D}^J_1$, $G^J_1$-invariant to the
action \eqref{dg}, \eqref{xxx}, can be written as:
\begin{equation}\label{aab}
-\ii\omega_{k\mu}(z,w) =2k\frac{\dd w \wedge \dd\bar{w}}{P^2} +
\mu\frac{\mc{A}\wedge \bar{\mc{A}}}{P},
~\mc{A}=\dd z+\dd w \bar{\eta}(z,w), ~P=1-w\bar{w},
\end{equation}
\begin{equation}\label{csv}
 ~z=\eta-w\bar{\eta},\quad\text{and~~~}~~~\eta = \eta(z,w):= \frac{z+\bar{z}w}{P}.
\end{equation}
$k$ indexes the positive discrete series of $\text{SU}(1,1)$
($2k\in\N$), while $\mu>0$ indexes the representations of the
Heisenberg group. 

The matrix of the balanced metric  \eqref{herm} $h=h(\varsigma)$, 
$\varsigma:=(z,w)\in\C\times\mc{D}_1$ is 
\begin{equation}\label{metrica}
h(\varsigma) =\left(\begin{array}{cc} \frac{\mu}{P} & \mu \frac{\eta}{P} \\
\mu\frac{\bar{\eta}}{P} &
\frac{2k}{P^2}+\mu\frac{|\eta|^2}{P}\end{array}\right).  
\end{equation}
The inverse of the matrix \eqref{metrica} reads
\begin{equation}\label{hinv}
h^{-1}(\varsigma)= \left(\begin{array}{cc}h^{z\bar{z}}&
      h^{z\bar{w}}\\h^{w\bar{z}}&h^{w\bar{w}}\end{array}\right)  =  \left(\begin{array}{cc}
    \frac{P}{\mu}+\frac{P^2|\eta|^2}{2k} & -\frac{P^2\eta}{2k} \\
-\frac{P^2\bar{\eta}}{2k} & \frac{P^2}{2k}\end{array}\right).
\end{equation}

 If we introduce the notation
\begin{equation}\label{DETG1}
\mc{G}_M(z):=\det
(h_{\alpha\bar{\beta}})_{\alpha,\beta=1,\dots,n},\end{equation}
then  we
find 
\begin{equation}\label{g311}
\mc{G}_{\mc{D}^J_1}(z,w)=\frac{2k\mu}{(1-w\bar{w})^3},~~z\in\C,~|w|<1.
\end{equation}

In the variables $(z,w)\in(\C,\mc{D}_1)$ the equations of geodesics
\eqref{GEOO}  read
\begin{equation}\label{geomic}
 \left\{
 \begin{array}{l}
 \frac{\dd^2 z}{\dd t^2}+\Gamma^z_{zz}\left(\frac{\dd z}{\dd
     t}\right)^2 +
2\Gamma^z_{zw}\frac{\dd z}{\dd t}  \frac{\dd w}{\dd t} +\Gamma
^z_{ww}\left(\frac{\dd w}{\dd t}\right)^2 =0   ;\\
  \frac{\dd^2 w}{\dd t^2}+\Gamma^w_{zz}\left(\frac{\dd z}{\dd
     t}\right)^2 +
2\Gamma^w_{zw}\frac{\dd z}{\dd t}  \frac{\dd w}{\dd t} +\Gamma
^w_{ww}\left(\frac{\dd w}{\dd t}\right)^2=0 .
     \end{array}
 \right.
\end{equation}
The equations \eqref{CRISTU} which determine the $\Gamma$-symbols for
the Siegel-Jacobi disk are 
\begin{equation}\label{ec528}
 \left\{
 \begin{array}{l}
h_{z\bar{z}}\Gamma^z_{zz}+h_{w\bar{z}}\Gamma^w_{zz}=
\frac{\pa h_{z\bar{z}}}{\pa z}; \\
h_{z\bar{w}}\Gamma^z_{zz}+h_{w\bar{w}}\Gamma^w_{zz}=
\frac{\pa h_{z\bar{w}}}{\pa z}. 
\end{array}
 \right.
\end{equation}

\begin{equation}\label{ec529}
 \left\{
 \begin{array}{l}
h_{z\bar{z}}\Gamma^z_{zw}+h_{w\bar{z}}\Gamma^w_{zw}=
\frac{\pa h_{w \bar{z}}}{\pa z}; \\
h_{z\bar{w}}\Gamma^z_{zw}+h_{w\bar{w}}\Gamma^w_{zw}=
\frac{\pa h_{w\bar{w}}}{\pa z}. 
\end{array}
 \right.
\end{equation}
\begin{equation}\label{ec530}
\left\{
 \begin{array}{l}
h_{z\bar{z}}\Gamma^z_{ww}+h_{w\bar{z}}\Gamma^w_{ww}=
\frac{\pa h_{w \bar{z}}}{\pa w}; \\
h_{z\bar{w}}\Gamma^z_{ww}+h_{w\bar{w}}\Gamma^w_{ww}=
\frac{\pa h_{w\bar{w}}}{\pa w}. 
\end{array}
 \right.
\end{equation}
From \eqref{ec528}-\eqref{ec530} or from \eqref{stefaN}, we get  
\begin{equation}\label{pdgM}
\begin{split}
\Gamma^z_{zz} & = h^{\bar{z}z}\frac{\pa h_{z\bar{z}}}{\pa z} +
h^{\bar{w}z}\frac{\pa h_{z\bar{w}}}{\pa z};~~
\Gamma^w_{zz}  =h^{\bar{z}w}\frac{\pa h_{z\bar{z}}}{\pa z}
+h^{\bar{w}w}\frac{\pa h_{z\bar{w}}}{\pa z};\\ 
\Gamma^z_{zw} &=h^{z\bar{z}}\frac{\pa
  h_{w\bar{z}}}{\pa z}+h^{w\bar{z}} \frac{\pa h_{w\bar{w}}}{\pa z};~~
\Gamma^w_{zw}  =h^{w\bar{w}}\frac{\pa
  h_{w\bar{w}}}{\pa z}+h^{z\bar{w}} \frac{\pa h_{w\bar{z}}}{\pa z};\\
\Gamma^z_{ww} & =h^{z\bar{z}}\frac{\pa
  h_{w\bar{z}}}{\pa w}+h^{w\bar{z}} \frac{\pa h_{w\bar{w}}}{\pa w};~~
\Gamma^w_{ww}  =h^{w\bar{w}}\frac{\pa
  h_{w\bar{w}}}{\pa w}+h^{z\bar{w}} \frac{\pa h_{w\bar{z}}}{\pa w}.
\end{split}
\end{equation}
With \eqref{metrica}, we calculate  easily the derivatives
\begin{equation}\label{ec531}
\begin{split}
\frac{\pa h_{z\bar{z}}}{\pa z } &= 0;~ \frac{\pa h_{z\bar{w}}}{\pa z}=
\frac{\mu}{P^2}; ~\frac{\pa h_{w\bar{z}}}{\pa z}=
\mu\frac{\bar{w}}{P^2}; \frac{\pa h_{w\bar{w}}}{\pa z }=\mu\frac{\bar{\eta}+\eta\bar{w}}{P^2};\\
\frac{\pa h_{w\bar{z}}}{\pa w } & = 2\mu\frac{\bar{w}\bar{\eta}}{P^2}; \frac{\pa h_{w\bar{w}}}{\pa w}=
\mu\frac{\bar{z}\bar{\eta}}{P^2}+3\mu\frac{\bar{w}|\eta|^2}{P^2}+4k\frac{\bar{w}}{P^3}.
\end{split}
\end{equation}
Introducing \eqref{ec531} into \eqref{ec528}-\eqref{ec530}, we finds
for   the  Christoffel's symbols $\Gamma$-s  the expressions
\begin{equation}\label{GAMM}
\begin{split}
\Gamma^z_{zz}  & =-\lambda\bar{\eta};~\Gamma^w_{zz}=\lambda; ~
\Gamma^z_{zw}=-\lambda\bar{\eta}^2+\frac{\bar{w}}{P};\\
 \Gamma^w_{wz} & =\lambda\bar{\eta};~
 \Gamma^z_{ww}=-\lambda\bar{\eta}^3; ~ \Gamma^w_{ww} =
 \lambda\bar{\eta}^2+2\frac{\bar{w}}{P}, ~ \lambda = \frac{\mu}{2k}.
\end{split}
\end{equation}

\begin{Proposition}\label{GER}
The equations of geodesics on the Siegel-Jacobi corresponding to
metric defined by $\omega_{k\mu}$  \eqref{aab} are
\begin{subequations}\label{geo}
\begin{align}
\mu\bar{\eta}G^2_1 & = 2k G_3, ~ G_1=\frac{\dd z}{\dd
  t}+\bar{\eta}\frac{\dd w}{\dd t}, ~ G_3 = \frac{\dd ^2z}{\dd
  t^2}+2\frac{\bar{w}}{P}\frac{\dd z}{\dd t}\frac{\dd w}{\dd t}, ~ P=1-w\bar{w};  \label{geo1}\\
\mu G^2_1 & = -2k G_2, ~G_2=\frac{\dd^2w}{\dd
  t^2}+2\frac{\bar{w}}{P}(\frac{\dd w}{\dd t})^2 . \label{geo2}
\end{align}
\end{subequations}
\end{Proposition}
The connection matrix $w_{\mc{D}^J_1} $ on $\mc{D}^J_1$ 
$$\theta_{\mc{D}^J_1}:=\left(\begin{array}{cc} \theta^z_z & \theta^z_w\\ \theta^w_z &
                                                                  \theta^w_w\end{array}\right)
                                                              =\left(\begin{array}{cc}\Gamma^z_{zz}\dd
                                                                       z
                                                                       +\Gamma^z_{zw}\dd
                                                                       w
                                                                    &
                                                                      \Gamma^z_{wz}+\Gamma^z_{ww}\dd
                                                                      w\\\Gamma^w_{zz}\dd
                                                                       z+\Gamma^w_{zw}\dd
                                                                       w&
                                                                          \Gamma^w_{wz}\dd
                                                                          z+\Gamma^w_{ww}\dd w\end{array}\right)$$
has the value
\begin{equation}\theta_{\mc{D}^J_1}=\left(\begin{array}{cc}
                          -\lambda\bar{\eta}\mc{A}+\frac{\bar{w}}{P}\dd
                          w & -\lambda\bar{\eta}^2\mc{A} +\frac{\bar{w}}{P}\dd
                          z\\ \lambda \mc{A} &
                          \lambda\bar{\eta}\mc{A}+2\frac{\bar{w}}{P}\dd
                                               w\end{array}\right) .\end{equation}
The covariant derivative of $\dd z$ on $\mc{D}^J_1$ has the expression 
\begin{equation}\label{DDZ}
\begin{split}
D (\dd z) & =\left(\begin{array}{cc}\dd z \dd w\end{array}\right) 
\left(\begin{array}{cc} \lambda\bar{\eta} &
                                            \lambda\bar{\eta}^2-\frac{\bar{w}}{P}\\
        \lambda\bar{\eta}^2-\frac{\bar{w}}{P} & \lambda \bar{\eta}^3\end{array}\right)
\left(\begin{array}{c}\dd z \\\dd w\end{array}\right)\\  & =\left(\begin{array}{cc}\mc{A} \dd w\end{array}\right) 
\left(\begin{array}{cc} \lambda\bar{\eta} &
                                            -\frac{\bar{w}}{P}\\
       -\frac{\bar{w}}{P} & 2\frac{\bar{\eta}\bar{w}}{P}\end{array}\right)
\left(\begin{array}{c}\mc{A} \\\dd w\end{array}\right).
\end{split}
\end{equation}
The covariante derivative of $\dd w$ has the expression
\begin{subequations}\label{DDW}
\begin{align}
- D (\dd w)  & =\left(\begin{array}{cc} \dd z \dd w\end{array}\right)
\left(\begin{array}{cc} \lambda & \lambda \bar{\eta} \\
        \lambda\bar{\eta} & \lambda\bar{\eta}^2+2\frac{\bar{w}}{P}\end{array}\right)
\left(\begin{array}{c}\dd z \\\dd w\end{array}\right)\\
& =  \left(\begin{array}{cc} \mc{A} \dd w\end{array}\right)
\left(\begin{array}{cc} \lambda & 0 \\
        0 &  2\frac{\bar{w}}{P}\end{array}\right)
\left(\begin{array}{c}\mc{A}\\\dd w\end{array}\right). 
\end{align} 
\end{subequations}

\subsection{Balanced metric on the Siegel-Jacobi ball}\label{BMSB}

The Jacobi algebra is the  the semi-direct sum
$\got{g}^J_n:= \got{h}_n\rtimes \got{sp}(n,\R )_{\C}$ \cite{JGSP,sbj,nou}.  The Heisenberg
algebra $\got{h}_n$ is generated  by    the boson creation 
(respectively, annihilation)
operators $\mb{a}^{\dagger}$, ($\mb{a}$), and 
\begin{equation}\label{baza1M}
[a_i,a^{\dagger}_j]=\delta_{ij}; ~ [a_i,a_j] = [a_i^{\dagger},a_j^{\dagger}]= 0 .
\end{equation}
  $K^{\pm,0}_{ij}$ are the
generators of the $\got{sp}(\R)_{\C}$ algebra 
\begin{subequations}\label{baza2M}
\begin{eqnarray}
 [K_{ij}^-,K_{kl}^-] & = & [K_{ij}^+,K_{kl}^+]=0 , ~2[K^-_{ij},K^0_{kl}]  =  K_{il}^-\delta_{kj}+K^-_{jl}\delta_{ki}\label{baza23}, \\
 2[K_{ij}^-,K_{kl}^+] & = & K^0_{kj}\delta_{li}+
K^0_{lj}\delta_{ki}+K^0_{ki}\delta_{lj}+K^0_{li}\delta_{kj}
\label{baza22}, \\
2[K^+_{ij},K^0_{kl}] & = & -K^+_{ik}\delta_{jl}-K^+_{jk}\delta_{li}
 \label{baza24}, ~
 2[K^0_{ji},K^0_{kl}] =  K^0_{jl}\delta_{ki}-K^0_{ki}\delta_{lj} .  
\end{eqnarray}
\end{subequations}
$\got{h}_n$ is an  ideal in $\got{g}^J_n$,  and 
\begin{subequations}\label{baza3M}
\begin{eqnarray}
\label{baza31}[a^{\dagger}_k,K^+_{ij}] & = & [a_k,K^-_{ij}]=0,  \\
~[a_i,K^+_{kj}]  & = &  
\frac{1}{2}\delta_{ik}a^{\dagger}_j+\frac{1}{2}\delta_{ij}a^{\dagger}_k ,~
 [K^-_{kj},a^{\dagger}_i]  = 
\frac{1}{2}\delta_{ik}a_j+\frac{1}{2}\delta_{ij}a_k , \\
~  [K^0_{ij},a^{\dagger}_k] & = & \frac{1}{2}\delta_{jk}a^{\dagger}_i,~
[a_k,K^0_{ij}]= \frac{1}{2}\delta_{ik}a_{j} .
\end{eqnarray}
\end{subequations}
 
Perelomov's CS vectors  \cite{perG}  associated to the group $G^J_n$ with 
Lie algebra the Jacobi algebra $\got{g}^J_n$ based the Siegel-Jacobi ball 
$ \mc{D}^J_n$
have been  defined as \cite{sbj,perG}
\begin{equation}\label{csuX}
e_{z,W}= \exp ({\mb{X}})e_0, 
~\mb{X} := \sqrt{\mu}\sum_{i=1}^n z_i \mb{a}^{\dagger}_i + \sum_{i,j =1}^n w_{ij}\mb{K}^+_{ij},~
 z\in \C^n;    W\in\mc{D}_n.
\end{equation} 

The vector $e_0$  appearing in (\ref{csuX})   verifies the relations 
\begin{equation}\label{vacuma}
\mb{a}_ie_o= 0,~ i=1, \cdots, n; ~
\mb{K}^+_{ij} e_0  \not=  0 ,~
\mb{K}^-_{ij} e_0 = 0 ,~
\mb{K}^0_{ij} e_0 =  \frac{k}{4}\delta_{ij} e_0, ~i,j=1,\dots,n .
\end{equation}
 The reproducing kernel $K(z,W)=(e_{z,W},e_{z,W})_{k\mu}$,
$z\in\C^n, W\in\mc{D}_n$ is 
\begin{equation}\label{kul}
K(z,W) = \det(M)~^{\frac{k}{2}}\exp\mu  F, M=(\un-W\bar{W})^{-1},
\end{equation}
\begin{subequations}
\begin{align}
2F & =2\bar{z}^tMz+z^t\bar{W}Mz+\bar{z}^tMW\bar{z},\label{FfF}\\
2F & = 2\bar{\eta}^t\eta -\eta^t\bar{W}\eta-\bar{\eta}^tW\bar{\eta},
\end{align}
\end{subequations}
\begin{equation}\label{etaZ}
\eta=M(z+W\bar{z}); ~ z=\eta-W\bar{\eta}.
\end{equation} 
The
manifold $\mc{D}^J_n$ has the K\"ahler potential \eqref{kelerX},
$f=\log K$, with $K$ given by \eqref{kul}, 
\begin{equation}\label{kelerX}
\begin{split}
f& \!=\!-\tfrac{k}{2}\log \det (\un-W\bar{W})\\
& +\mu \{
\bar{z}^t(\un-W\bar{W})^{-1}{z}
+ \tfrac{1}{2}z^t[\bar{W}(\un-W\bar{W})^{-1}]z
+\tfrac{1}{2}\bar{z}^t [ (\un-W\bar{W})^{-1}W]\bar{z} \}.
\end{split}
\end{equation}
We use the following notation for the matrix of the metric:
\begin{equation}\label{math}
\!h=\!\left(\begin{array}{cc} h_1 &h_2\\h_3& h_4\end{array}\right)\!=\!
\left(\begin{array}{cc} h_{i\bar{j}}& h_{i\bar{p}\bar{q}} \\
    h_{pq\bar{i}} &  h_{pq\bar{m}\bar{n}} \end{array}\right)\!\in
M(n(n+3)/2,\C),~p\leq q, m\leq n.
 ~ h=h^{*},\end{equation}

We use the following convention: indices of $z\in M(n\times 1,\C)$ are
denoted with: $i,j,k,l$; indices of
$W=W^t, W\in M(n,\C)$ are denoted with:  $m,n,p,q,r,s,t,u,v$.

We have
determined the ``inverse''  $h^{-1}$ of \eqref{math} such that
\begin{equation}\label{ciudat}
\left(\begin{array}{cc} h_1 &h_2\\h_3& h_4\end{array}\right)
\left(\begin{array}{cc}h^1&h^2\\h^3&h^4\end{array}\right)=
\left(\begin{array}{cc} \un & 0\\ 0 & \Delta \end{array}\right),
\end{equation}
where $\Delta$ is defined in \eqref{mircea}
(see equation (4.5) in \cite{nou}):
\begin{equation}\label{mircea}
\Delta^{ij}_{pq}:=\frac{\pa w_{ij}}{\pa w_{pq}} =
\delta_{ip}\delta_{jq}+\delta_{iq}\delta_{jp}-\delta_{ij}\delta_{pq}\delta_{ip}=
(\delta_{ip}\delta_{jq}+\delta_{iq}\delta_{jp})f_{ij},
~ w_{ij}=w_{ji} .
\end{equation}
We use  the notation
\begin{equation}\label{fqq}
 f_{pq}:=1-\frac{1}{2}\delta_{pq}; ~~~d_{pq}:=1-\delta_{pq}.
\end{equation}

In \cite{SB15} we have proved
\begin{Theorem}\label{mainTH}
The  K\"ahler two-form $\omega_{\mc{D}^J_n}$, associated with
 the  balanced metric  of the Siegel-Jacobi ball   $\mc{D}^J_n$,  $G^J_n$-invariant to the action
\eqref{x44}, \eqref{xxxX}  has the expression
\begin{equation}\label {aabX}
\begin{split}
- \ii\omega_{\mc{D}^J_n}(z,W)&\! =\! \tfrac{k}{2}\tr (\mc{B}\wedge\bar{\mc{B}})\!
 +\mu \tr (\mc{A}^t\bar{M}\wedge \bar{\mc{A}}), ~\mc{A} =\dd z+ \dd W\bar{\eta},\\
\mc{B} & = M\dd W,~    M  =\!(\un-W\bar{W})^{-1}.
\end{split}
\end{equation}
The matrix \eqref{math} of the metric on $\mc{D}^J_n$  has the matrix elements
\eqref{hcomp}:
\begin{subequations}\label{hcomp}
\begin{align}
h_{i\bar{j}} & =\mu \bar{M}_{ij},\label{hc11}\\
h_{i\bar{p}\bar{q}} & = \mu (\eta_q\bar{M}_{ip}
+\eta_p\bar{M}_{iq})f_{pq},\label{hc12}\\
h_{pq\bar{i}} & = \mu (\bar{\eta}_q\bar{M}_{pi} +\bar{\eta}_p\bar{M}_{qi})f_{pq},\label{hc21}\\
h_{pq\bar{m}\bar{n}} & =\frac{k}{2}h^k_{pq\bar{m}\bar{n}}+\mu
h^{\mu}_{pq\bar{m}\bar{n}},\label{hc22}\\
h^k_{pq\bar{m}\bar{n}} & =
2(M_{mp}M_{nq}+M_{mq}M_{np})-2M_{mp}(M_{np}\delta_{pq}+M_{mq}\delta_{mn})\label{hc23}
 + M^2_{mp}\delta_{pq}\delta_{mn} \\
&= 2M_{mp}M_{nq}d_{pq}+2M_{mq}M_{np}d_{mn}+  
M^2_{mp}\delta_{pq}\delta_{mn} ,\nonumber\\
h^{\mu}_{pq\bar{m}\bar{n}} & = \bar{\eta}_p(\eta_n M_{mq}+\eta_m
M_{nq})d_{pq}+\bar{\eta}_q({\eta}_nM_{mp}+\eta_mM_{np}) \label{hc24}\\
& ~~ -{\eta}_m(\bar{\eta}_pM_{mq}+\bar{\eta}_qM_{mp})\delta_{mn}
+\bar{\eta}_p\eta_mM_{mp}\delta_{pq}\delta_{mn}\nonumber\\
& = [\bar{\eta}_p(\eta_n\bar{M}_{qm}+\eta_m\bar{M}_{qn})+\bar{\eta}_q(\eta_n\bar{M}_{pm}+\eta_m\bar{M}_{pn})]f_{pq}f_{mn}.
\nonumber
\end{align}
\end{subequations}    
The ``inverse''     of the matrix with elements $h^k_{pq\bar{m}\bar{n}}$ given by \eqref{hc23}
has the matrix elements 
\begin{equation}\label{kpqmn}
  k_{mn\bar{u}\bar{v}}=\frac{1}{2}(N_{vn}\bar{N}_{mu}+N_{vm}\bar{N}_{nu}),
  ~N=\un-W\bar{W},\end{equation}
and we have 
\begin{equation}\label{ceinv}
\sum_{m\le n}h^k_{pq\bar{m}\bar{n}}k_{mn\bar{u}\bar{v}}=\Delta^{uv}_{pq}.
\end{equation}
The ``inverse''  $h^{-1}$ of the metric matrix $h$
  which verifies \eqref{ciudat},
 has the elements $h^1$-- $h^4$
given by 
\begin{subequations}\label{nr11}
\begin{align}
(h^1)_{ij} &
             =\sigma\bar{M}^{-1}_{ij},
             \quad \sigma = \frac{1}{\mu}+\alpha\frac{2}{k},~
\alpha=\eta^t\bar{M}^{-1}\bar{\eta}=\bar{\eta}_nS_n,~~S_n=\sum\eta_q\bar{M}^{-1}_{qn};\label{NR1}\\
(h^2)_{i\bar{m}\bar{n}} &
                         =-\frac{1}{k}(S_n\bar{M}^{-1}_{im}+S_m\bar{M}^{-1}_{in});\label{NR2}\\
(h^3)_{mn\bar{i}} &
                    =-\frac{1}{k}(\bar{S}_n\bar{M}^{-1}_{mi}+\bar{S}_m\bar{M}^{-1}_{ni});\label{NR3}\\
(h^4)_{pq\bar{m}\bar{n}} & =
\frac{2}{k}(h^k)^{-1}_{pq\bar{m}\bar{n}}=\frac{1}{k}(\bar{M}^{-1}_{qn}\bar{M}^{-1}_{pm}+\bar{M}^{-1}_{pn}\bar{M}^{-1}_{qm}). \label{NR4}
\end{align}
\end{subequations}
The determinant of the metric  matrix $h$ is
\begin{equation}\label{DTH}
\mc{G}_{\mc{D}^J_n}(z,W)= \det h _{\mc{D}^J_n}(z,W)= (\frac{k}{2})^{\frac{n(n+1)}{2}}\mu^n\det
(\un-W\bar{W})^{-(n+2)}.\end{equation}
In the case $n=1$ formulae \eqref{hcomp}, \eqref{nr11}, \eqref{DTH}
became the formulae \eqref{metrica}, \eqref{hinv}, respectively
\eqref{g311}, with $2k\leftrightarrow\frac{k}{2}$. 
\end{Theorem}
Now we recall  (see Proposition 4 in \cite{nou}, Lemma 2 and
Proposition 3 in \cite{csg}, and Propositions 1 and 2 in\cite{SB14})  the significance of the change of coordinates
\eqref{etaZ} as a $FC$-transform  in the language of  coherent states
and also in the context of fundamental conjecture for homogeneous
K\"ahler manifolds of Gindikin-Vinberg \cite{GV} (see the proof in
\cite{DN}).
\begin{Proposition}\label{expl}
a) The change of coordinates \eqref{etaZ} $(z,W)=FC(\eta,W)$,
$z=\eta-W\bar{\eta}$ is a $FC$-transform in the meaning of \eqref{etild}.

b) The $FC$-transform \eqref{etaZ}
$\C^n\times\mc{D}_n\ni(\eta,W)\stackrel{\mbox{FC}}{\rightarrow}(z,W)\in\mc{D}^J_n$
is  K\"ahler homogeneous 
diffeomorphism, i.e.
$$\omega_{\C^n\times\mc{D}_n}(z,\eta)=FC^*[\omega_{\mc{D}^J_n}(z,W)]=\omega_{\mc{D}_n}(W)+\omega_{\C^n}(\eta),$$
where 
\begin{equation}\label{2omega}
-\ii\omega_{\mc{D}_n}(W)=\frac{k}{2}\tr (\mc{B}\wedge\bar{\mc{B}}),~~ -\ii
\omega_{\C^n}=\mu\tr(\dd \eta^t\wedge\dd \bar{\eta}).\end{equation}
The K\"ahler two-form $\omega_{\C^n\times\mc{D}_n}(z,\eta)$ is
invariant to the $G^J_n$-action on $\mc{D}_n\times\C^n$:
$(g,\alpha)\times (\eta,W)\rightarrow(\eta_1,W_1)$, where $g$ has the
expression \eqref{dgM}, $W_1$ is given by \eqref{x44}, and 
$$\eta_1=p(\eta+\alpha)+q(\bar{\eta}+\bar{\alpha}).$$
\end{Proposition}
In the following proposition we collect several  geometric properties
of the Siegel-Jacobi group and its action on the Siegel-Jacobi ball
(see Remark 4, Proposition 3 and Proposition 4 from \cite{SB15}, where
all the terms appearing in the enunciation below are explained).
\begin{Proposition}\label{TOTU}
i) The Jacobi group $G^J_n$  is an unimodular, non-reductive, algebraic group of
Harish-Chandra type. \\
ii) The Siegel-Jacobi domain $\mc{D}^J_n$  is a  homogeneous reductive,
non-symmetric manifold associated to the Jacobi group $G^J_n$ by the
generalized Harish-Chandra embedding.\\
iii) The homogeneous K\"ahler manifold $\mc{D}^J_n$ is contractible. \\
iv) The  K\"ahler potential  of the Siegel-Jacobi ball 
is global. $\mc{D}^J_n$ is  a Lu Qi-Keng  manifold, with
 nowhere vanishing diastasis.\\
v) The manifold $\mc{D}^J_n$ is a quantizable manifold.\\
vi) The $\mc{D}^J_n$ is projectively induced, and the Jacobi group $G^J_n$
is a CS-type group.\\
vii) The Siegel-Jacobi ball $\mc{D}^J_n$ is not an
  Einstein metric with respect to the balanced  metric
  attached to the K\"ahler two-form \eqref{aabX}, but it is one with respect to the Bergman metric
  corresponding to the Bergman K\"ahler two-form. \\
ix) The scalar curvature is constant and negative. 
\end{Proposition}

\section{Geodesics on the Siegel-Jacobi ball}\label{P2}
\subsection{Calculation of $\Gamma$-s}\label{S1}
In the case of the Siegel-Jacobi ball, the equations of geodesics
\eqref{GEOO} read
\begin{subequations}\label{2EE}
\begin{align}
\frac{\dd^2z_i}{\dd t^2}+\Gamma^{i}_{jk}\frac{\dd z_j}{\dd
  t}\frac{\dd z_k}{\dd t}+2\sum_{p\le q}\Gamma^i_{jpq}\frac{\dd z_j}{\dd
  t}\frac{\dd w_{pq}}{\dd t}+\sum_{p\le q, m\le n}\Gamma^i_{pqmn}\frac{\dd w_{pq}}{\dd
  t}\frac{\dd w_{mn}}{\dd t}   & = 0,\label{2EE1}\\
\frac{\dd^2 w_{pq}}{\dd t^2}+\Gamma^{pq}_{jk}\frac{\dd z_j}{\dd
  t}\frac{\dd z_k}{\dd t}+2\sum_{m\le n}\Gamma^{pq}_{imn}\frac{\dd z_i}{\dd
  t}\frac{\dd w_{mn}}{\dd t}
+\sum_{m\le n, u\le v}\Gamma^{pq}_{mnuv}\frac{\dd w_{mn}}{\dd t}\frac{\dd w_{uv}}{\dd t}\label{2EE2} &
=0 .
\end{align}
\end{subequations}
Equations \eqref{geomic} are a particular case of \eqref{2EE}.

The equations \eqref{CRISTU} which determine the $\Gamma$-symbols for
the Siegel-Jacobi ball $\mc{D}^J_n$  are 
\begin{equation}\label{ec5288}
 \left\{
 \begin{array}{l}
h_{i\bar{j}}\Gamma^i_{kl}+\sum_{m\le n}h_{mn\bar{j}}\Gamma^{mn}_{kl}=
\frac{\pa h_{l\bar{j}}}{\pa z_k}=\frac{\pa h_{k\bar{j}}}{\pa z_l} ; \\
h_{i\bar{p}\bar{q}}\Gamma^{i}_{kl}+\sum_{m\le n}h_{mn\bar{p}\bar{q}}\Gamma^{mn}_{kl}=
\frac{\pa h_{k\bar{p}\bar{q}}}{\pa z_l}. 
\end{array}
 \right.
\end{equation}
\begin{equation}\label{ec5298}
 \left\{
 \begin{array}{l}
h_{i\bar{j}}\Gamma^i_{kpq}+\sum_{m\le n}h_{mn\bar{j}}\Gamma^{mn}_{kpq}=
\frac{\pa h_{pq\bar{j}}}{\pa z_k} ; \\
h_{i\bar{p}\bar{q}}\Gamma^{i}_{krs} +\sum_{m\le n}h_{mn\bar{j}}\Gamma^{mn}_{kpq}=
\frac{\pa h_{pq\bar{j}}}{\pa z_k}.\\
\end{array}
 \right.
\end{equation}
\begin{equation}\label{ec5299}
 \left\{
 \begin{array}{l}
h_{i\bar{j}}\Gamma^i_{pquv}+\sum_{m\le n}h_{mn\bar{j}}\Gamma^{mn}_{pquv}=
\frac{\pa h_{uv\bar{j}}}{\pa w_{pq}} ; \\
h_{i\bar{s}\bar{t}}\Gamma^{i}_{pquv}+\sum_{m\le n}h_{mn\bar{s}\bar{t}}\Gamma^{mn}_{pquv}=
\frac{\pa h_{uv\bar{s}\bar{t}}}{\pa w_{pq}}. 
\end{array}
 \right.
\end{equation}
Equations \eqref{ec5288}-\eqref{ec5299} are a generalization on
$\mc{D}^J_n$ of the
corresponding ones \eqref{ec528}-\eqref{ec530} on $\mc{D}^J_1$.

With formula \eqref{stefaN}, we  have 
\begin{subequations}\label{pdG}
\begin{align}
\Gamma^i_{jk}
 & = h^{l\bar{i}}\frac{\pa h_{k\bar{l}}}{\pa z_j}+ \sum_{p\le q}h^{pq\bar{i}}\frac{\pa
  h_{k\bar{p}\bar{q}}}{\pa z_j},\\
\Gamma^{i}_{jpq} & = h^{k\bar{i}}\frac{\pa h_{pq\bar{k}}}{\pa
    z_j}+\sum_{r\le s}h^{rs\bar{i}}\frac{\pa h_{pq\bar{r}\bar{s}}}{\pa z_j}, \\
\Gamma^{pq}_{jk} & =h^{i\bar{p}\bar{q}}\frac{\pa h_{k\bar{i}}}{\pa
    z_j}+\sum_{m\le n}h^{mn\bar{p}\bar{q}}\frac{\pa h_{k\bar{m}\bar{n}}}{\pa z_j},
  \\
\Gamma^{pq}_{imn}  & =h^{j\bar{p}\bar{q}}\frac{\pa h_{mn\bar{j}}}{\pa
  z_i}+\sum_{u\le v}h^{uv\bar{p}\bar{q}}\frac{\pa h_{mn\bar{u}\bar{v}}}{\pa z_i},
\\\Gamma^i_{pqmn}  & =h^{j\bar{i}}\frac{\pa h_{mn\bar{j}}}{\pa
  w_{pq}}+\sum_{r\le s}h^{rs\bar{i}}\frac{\pa h_{mn\bar{r}\bar{s}}}{\pa w_{pq}}, \label{724e}\\
\Gamma^{pq}_{mnuv}  & =h^{i\bar{p}\bar{q}}\frac{\pa h_{uv\bar{i}}}{\pa
  w_{mn}}+\sum_{s\le t}h^{st\bar{p}\bar{q}}\frac{\pa h_{uv\bar{s}\bar{t}}}{\pa
  w_{mn}}.\label{724f} 
\end{align}
\end{subequations}
Equations \eqref{pdG} for $\mc{D}^J_n$ generalize equations
\eqref{pdgM} for $\mc{D}^J_1$.

We shall use the formulae given in Theorem \ref{mainTH} for the 
matrix elements of the metric and its inverse and also some formulae  (see (3.28) and (3.29)
in \cite{SB15})
\begin{subequations}\label{LIKK}
\begin{align}
\frac{\pa M_{ab} }{\pa w_{ik}} & =  (M_{ai}X_{kb}+
M_{ak}X_{ib})f_{ik},~~\text{where}~~X=X^t=\bar{W}M=\bar{M}\bar{W},\label{LIK}\\
\frac{\pa X_{ab} }{\pa w_{ik}}  & = 
(X_{ai}X_{bk}+X_{ak}X_{ib})f_{ik} , \\
\frac{\pa \bar{ X}_{ab} }{\pa w_{ik}}  &  = 
(M_{ai}M_{bk}+M_{ak}M_{bi})f_{ik} \label{LIK2}.
\end{align}
\end{subequations}
\begin{equation}\label{morD}\begin{split}
\frac{\pa \eta_q}{\pa z_l} & = M_{ql},~~
\frac{\pa \bar{\eta}_q}{\pa z_j}  ={X}_{qj}, \\
\frac{\pa \eta_t}{\pa w_{pq}} & =(M_{tp}\bar{\eta}_q +M_{tq}\bar{\eta}_p)f_{pq}, ~~
\frac{\pa \bar{\eta}_n}{\pa w_{pq}}  =
(\bar{\eta}_pX_{qn}+\bar{\eta}_qX_{pn})f_{pq} .
\end{split}
\end{equation}

Now we want to calculate the partial derivatives which appear in
\eqref{pdG}, a generalization of the partial derivatives \eqref{ec531}.  We use \eqref{LIKK} and \eqref{morD} to calculate the partial
 derivatives which appear in the equations \eqref{pdG} of $\Gamma$-s
 and we get
\begin{subequations}
\begin{align*}
\frac{\pa h_{k\bar{l}}}{\pa z_j} & = 0,~~
 \frac{\pa  h_{k\bar{p}\bar{q}}}{\pa z_j} = \mu
                                            (M_{pk}M_{qj}+M_{pj}M_{qk})f_{pq},\\
\frac{h_{pq\bar{k}}}{\pa z_j} &= \mu(X_{qj}M_{kp}+X_{pj}M_{kq}) f_{pq},~~
\frac{\pa h_{pq\bar{r}\bar{s}}}{\pa z_j}  = \mu
  M_{rp}[{X}_{qj}\eta_s+\bar{\eta}_qM_{sj}],\\
\!\frac{1}{\mu}\frac{\pa h_{pq\bar{r}\bar{s}}}{\pa z_j} &\!=\!
                                                        [X_{pj}(\eta_s\bar{M}_{qr}+\eta_{r}\bar{M}_{qs})+X_{qj}(\eta_s\bar{M}_{pr}+\eta_{r}\bar{M}_{ps})]f_{pq}f_{rs}\\
& \!+[\bar{\eta}_p (M_{sj}\bar{M}_{qr}+M_{rj}\bar{M}_{qs})+\bar{\eta}_q(M_{sj}\bar{M}_{pr}+M_{rj}\bar{M}_{ps})]f_{pq}f_{rs},\\
\!\frac{1}{\mu}\frac{\pa h_{mn\bar{j}}}{\pa
  w_{pq}} &
\!=\! [\Lambda^{pq}_n\bar{M}_{mj} +
   \bar{\eta}_n(M_{jp}X_{qm}+M_{jq}X_{pm})]f_{mn}f_{pq}\\
\!& \!+\! [ \Lambda^{pq}_m\bar{M}_{nj} +
   \bar{\eta}_m(M_{jp}X_{qn}+M_{jq}X_{pn})]f_{mn}f_{pq},\\
\!\frac{\pa h^k_{mn\bar{r}\bar{s}}}{\pa w_{pq}} &
\!=\!\!2[\!(M_{rp}X_{qm}\!+\!M_{rq}X_{pm}\!)M_{sn}\!+\!M_{rm}(M_{sp}X_{qn}\!+\!M_{sq}X_{pn})\!]d_{mn}f_{pq}\\
\!&\!+\!2[\!(M_{rp}X_{qn}\!+\!M_{rq}X_{pn})M_{sm}\!+\!M_{rn}(M_{sp}X_{qm}\!+\!M_{sq}X_{pm})]d_{rs}f_{pq}\\
\!&\!+\!2M_{rn}(M_{sp}X_{qm}+M_{sq}X_{pm})\delta_{mn}\delta_{rs}f_{pq},
\end{align*}
where
\end{subequations}
\begin{equation}\label{LPQ}
\Lambda^{pq}_m:=\bar{\eta}_pX_{qm}+\bar{\eta}_qX_{pm}.
\end{equation}
We have also
\begin{equation}
\frac{\pa h^{\mu}_{mn\bar{r}\bar{s}}}{\pa w_{pq}}=\Upsilon^{pq}_{mn}f_{mn}f_{rs}f_{pq},
\end{equation}
where
\begin{subequations}
\begin{align*}
\Upsilon^{pq}_{mn}&= \Lambda^{pq}_m(\eta_s\bar{M}_{nr}+\eta_r\bar{M}_{ns}) 
+\Lambda^{pq}_n(\eta_s\bar{M}_{mr}+\eta_r\bar{M}_{ms}) \\
&+\bar{\eta}_m[(\bar{\eta}_qM_{sp}+\bar{\eta}_pM_{sq})\bar{M}_{nr}
  +\eta_s(M_{rp}X_{qn}+M_{rq}X_{pn})]\\&+\bar{\eta}_m[(\bar{\eta}_qM_{rp}+\bar{\eta}_pM_{rq})\bar{M}_{ns}+\eta_r(M_{sp}X_{qn}+M_{sq}X_{pn})]\\
&+\bar{\eta}_n[(\bar{\eta}_qM_{sp}+\bar{\eta}_pM_{sq})\bar{M}_{mr}
  +\eta_s(M_{rp}X_{qm}+M_{rq}X_{pm})]\\
&+\bar{\eta}_n[(\bar{\eta}_qM_{rp}+\bar{\eta}_pM_{rq})\bar{M}_{ms}+\eta_r(M_{sp}X_{qm}+M_{sq}X_{pm})].
\end{align*}
\end{subequations}

In order to calculate $\Gamma^i_{pqmn}$ with \eqref{724e}, we
calculate the first term 
 $$
_1\Gamma^i_{pqmn} :=h^{j\bar{i}}\frac{\pa h_{mn\bar{j}}}{\pa
  w_{pq}}=(1+2\epsilon\alpha)
X^i_{pqmn}f_{mn},~~\epsilon=\frac{\mu}{k},$$
where
\begin{equation}
X^i_{pqmn}  =(\delta_{mi}\Lambda^{pq}_m+\delta_{ni}\Lambda^{pq}_n +\bar{\eta}_m \Omega^{i}_{pqn}
  +\bar{\eta}_n\Omega^{i}_{pqm})f_{pq},
\end{equation}
\begin{equation}\label{OPQ}
\Omega^{i}_{pqm}:=\delta_{ip}X_{qm}+\delta_{iq}X_{pm}.
\end{equation}
Now we calculate
$$_2\Gamma^i_{pqmn}:=\sum_{r\le s}h^{rs\bar{i}}\frac{\pa
  h_{mn\bar{r}\bar{s}}}{\pa
  w_{pq}}=_{2k}\Gamma^i_{pqmn}+_{2\mu}\Gamma^i_{pqmn},$$
where
$$_{2k}\Gamma^i_{pqmn}= -\sum_{s\le r}(\bar{S}_s\bar{M}^{-1}_{ri}+\bar{S}_r\bar{M}^{-1}_{si} )\frac{\pa
  h^k_{mn\bar{r}\bar{s}}}{\pa w_{pq}}.$$
We
find $$_{2k}\Gamma^i_{pqnn}=-\frac{1}{2}X^i_{pqnn},$$
while  $$_{2k}\Gamma^i_{pqmn} =-X^i_{pqmn},~ m\not=n,$$
i.e.
$$_{2k}\Gamma^i_{pqmn}=-X^i_{pqmn}f_{mn}.$$
For $$_{2\mu}\Gamma^i_{pqmn}:=\sum_{r\le s}h^{st\bar{p}\bar{q}}\frac{\pa h^{\mu}_{uv\bar{s}\bar{t}}}{\pa
  w_{mn}},$$
we find $$_{2\mu}\Gamma^i_{pqmn}:=-\epsilon f_{mn}f_{pq}
Y^i_{pqmn},$$
where 
\begin{equation}
\begin{split}
Y^i_{pqmn} &
=\alpha(\delta_{ni}\Lambda^{pq}_m+\delta_{mi}\Lambda^{pq}_n+\bar{\eta}_m\Omega^i_{pqn}+\bar{\eta}_n\Omega^i_{pqm})
\\
& +
2S_i(\bar{\eta}_n\Lambda^{pq}_m+\bar{\eta}_m\Lambda^{pq}_m)+2\bar{\eta}_p\bar{\eta}_q(\bar{\eta}_m\delta_{in}+\bar{\eta}_n\delta_{im})\\
& +2\bar{\eta}_n\bar{\eta}_m(\bar{\eta}_q\delta_{ip}+\bar{\eta}_p\delta_{iq} ).
\end{split}
\end{equation}
We find 
$$\Gamma^i_{pqmn}=\epsilon f_{mn} f_{pq}Z^i_{pqmn},$$
where
\begin{equation}\label{ZI}
\begin{split}
Z^i_{pqmn} &=
\alpha(\delta_{ni}\Lambda^{pq}_m+\delta_{mi}\Lambda^{pq}_n+\bar{\eta}_m\Omega^i_{pqn}+\bar{\eta}_n\Omega^i_{pqm})\\
& - 2S_i(\bar{\eta}_n\Lambda^{pq}_m+\bar{\eta}_m\Lambda^{pq}_n)-2\bar{\eta}_p\bar{\eta}_q(\bar{\eta}_m\delta_{in}+\bar{\eta}_n\delta_{im})\\
& - 2\bar{\eta}_n\bar{\eta}_m(\bar{\eta}_q\delta_{ip}+\bar{\eta}_p\delta_{iq}).
\end{split}
\end{equation}
In \eqref{724f} we calculate firstly
\begin{equation}
\begin{split}
_1 \Gamma^{pq}_{mnuv}  & :=h^{i\bar{p}\bar{q}}\frac{\pa h_{uv\bar{i}}}{\pa
  w_{mn}}\\
& = -\epsilon f_{uv}f_{mn}S_q(\delta_{pu}\Lambda^{mn}_v+\bar{\eta}_v\Omega^p_{mnu}+\delta_{pv}\Lambda^{mn}_u+\bar{\eta}_u\Omega^p_{mnv})\\
&  -\epsilon f_{uv}f_{mn}S_q S_p(\delta_{qu}\Lambda^{mn}_v+\bar{\eta}_v\Omega^q_{mnu}+\delta_{qv}\Lambda^{mn}_u+\bar{\eta}_u\Omega^q_{mnv}).
\end{split}
\end{equation}
Now we calculate in \eqref{724f}
\begin{equation}
_{2k}\Gamma^{pq}_{mnuv} :=h^{rs\bar{i}}\frac{\pa
  h^k_{mn\bar{r}\bar{s}}}{\pa w_{pq}} = f_{uv}f_{mn}(\delta_{vq}\Omega^p_{mnv}+\delta_{up}\Omega^q_{mnv}+\delta_{uq}\Omega^p_{mnv}+\delta_{vp}\Omega^q_{mnu}).
\end{equation}

$$_{2\mu}\Gamma^{pq}_{mnuv}:=   h^{st\bar{p}\bar{q}}\frac{\pa h^{\mu}_{uv\bar{s}\bar{t}}}{\pa
  w_{mn}}=\epsilon f_{uv}f_{mn}U^{st}_{mnuv},$$
where 
\begin{equation}
\begin{split}
U^{st}_{mnuv} & = \Lambda^{mn}_u(S_q\delta_{pv}+S_{p}\delta_{qv})+
\Lambda^{mn}_v(S_q\delta_{pu}+S_{p}\delta_{qu})\\
 &~~+
 (\bar{\eta}_n\delta_{mq}+\bar{\eta}_m\delta_{nq})+(\delta_{vp}\bar{\eta}_u+\delta_{up}\bar{\eta}_v)\\
&~~+
(\bar{\eta}_n\delta_{mp}+\bar{\eta}_m\delta_{np})(\delta_{vq}\bar{\eta}_u+\delta_{uq}\bar{\eta}_v)\\
&~~+ \bar{\eta}_u(S_q\Omega^p_{mnv}+S_p\Omega^q_{mnv})+ \bar{\eta}_v(S_q\Omega^p_{mnu}+S_p\Omega^q_{mnu}).
\end{split}
\end{equation}
Now we calculate
$$\Gamma^{pq}_{mnuv}:=_{2k}\Gamma^{pq}_{mnuv}+_{2\mu}\Gamma^{pq}_{mnuv}.$$
We find
$$\Gamma^{pq}_{mnuv}=f_{uv}f_{mn}V^{pq}_{mnuv},$$
\begin{equation}\label{VPQ}
\begin{split}
V^{pq}_{mnuv} &
=\delta_{vq}\Omega^p_{mnu}+\delta_{up}\Omega^q_{mnv}+\delta_{uq}\Omega^p_{mnv}+\delta_{vp}\Omega^q_{mnu}+\epsilon
  W^{pq}_{mnuv},\\
W^{pq}_{mnuv} & =
(\bar{\eta}_n\delta_{mq}+\bar{\eta}_m\delta_{nq})(\delta_{vp}\bar{\eta}_u+\delta_{up}\bar{\eta}_{v})
 + (\bar{\eta}_n\delta_{mp}+\bar{\eta}_m\delta_{np})(\delta_{vq}\bar{\eta}_u+\delta_{uq}\bar{\eta}_{v}).
\end{split}
\end{equation}
We have proved
\begin{Proposition}\label{CRR}
The Christoffel's symbols of the Siegel-Jacobi ball $\mc{D}^J_n$
endowed with the balanced metric attached to the K\"ahler two-form \eqref{aabX}
  are 
\begin{subequations} \label{CRMARE}
\begin{align}
\Gamma^{i}_{jk} &
=-\epsilon(\bar{\eta}_j\delta_{ik}+\bar{\eta}_k\delta_{ij}),\quad\emph{\text{where~~}}
\epsilon=\frac{\mu}{k},\\
\Gamma^{i}_{jpq} &
=[(1+\epsilon\alpha)(X_{qj}\delta_{pi}+X_{pj}\delta_{qi}) -\epsilon
S_i(\bar{\eta}_pX_{qj}+\bar{\eta}_qX_{pj})]f_{pq}\\
 & -\epsilon[\bar{\eta}_j(\bar{\eta}_q\delta_{pi}+\bar{\eta}_p\delta_{qi})+2\bar{\eta}_p\bar{\eta}_q\delta_{ij}]f_{pq},\\
\Gamma^{pq}_{jk} & =\epsilon (\delta_{kp}\delta_{jq}+\delta_{pj}\delta_{qk}), \\
\Gamma^{pq}_{imn} &
=\epsilon[\bar{\eta}_m(\delta_{np}\delta_{iq}+\delta_{nq}\delta_{ip})+\bar{\eta}_n(\delta_{mp}\delta_{iq}+\delta_{mq}\delta_{ip})]f_{mn},\\
\Gamma^i_{pqmn} & =\epsilon f_{mn}f_{pq} Z^i_{pqmn},\\
\Gamma^{pq}_{mnuv}& =f_{uv}f_{mn}V^{pq}_{mnuv}.
\end{align}
\end{subequations}
$S_i$ was defined in \eqref{NR1}, $X$ was defined in \eqref{LIK},  $Z^i_{pqmn}$ is given by \eqref{ZI}, $V^{pq}_{mnuv}$ is given by
\eqref{VPQ}, $\Lambda^{pq}_m$ is defined in \eqref{LPQ},
$\Omega^{i}_{pqm}$ is given by \eqref{OPQ}.

Equations \eqref{CRMARE} of the Cristolffel's  symbols for
$\mc{D}^J_n$ generalize the corresponding ones \eqref{GAMM} on $\mc{D}^J_1$.
\end{Proposition}
For the proof of the last assertion in Proposition \ref{CRR}, we have
in the case case $n=1$ the following relations:\\
$\alpha= |\eta|^2P$, $S=\eta P$, $M=\bar{M}=P^{-1}$,
$X:=\bar{W}M\rightarrow\bar{w}P^{-1}$, $\Lambda=\bar{\eta}\bar{w}/P$,
$\Omega=X$,
$\Upsilon=\frac{4\bar{\eta}}{P^2}(\bar{\eta}+2\eta\bar{w})$,
  $X^i_{pqmn}\rightarrow 2(\Lambda+\bar{\eta}\Omega)=4\frac{\bar{\eta}\bar{w}}{P}$,
  $Y=4\bar{\eta}^2(\bar{\eta}+2\bar{w}\eta)$, $Z=-4\bar{\eta}^3$.\\
  $\frac{1}{2k}$ corresponds in the case $n=1$ to $\frac{2}{k}$, i.e. $\epsilon=\frac{1}{2}\lambda$. 

\subsection{Equations of geodesics}\label{S3}
Now we calculate the equations of geodesics \eqref{2EE}. In order to calculate
\eqref{2EE1}, we find successively:
$$\Gamma^{i}_{jk}\frac{\dd z_j}{\dd
  t}\frac{\dd z_k}{\dd t}= -2\epsilon L\frac{\dd z_i}{\dd t},$$
$$\sum_{p\le q}\Gamma^i_{jpq}\frac{\dd z_j}{\dd
  t}\frac{\dd w_{pq}}{\dd t} =
Q_i+\epsilon[\alpha Q_i-S_i\bar{\eta}^tQ-L\lambda_i-\bar{\eta}^t\frac{\dd W}{\dd t}\bar{\eta}\frac{\dd
  z_i}{\dd t}],$$
$$\sum_{p\le q, m\le n}\Gamma^i_{pqmn}\frac{\dd w_{pq}}{\dd
  t}\frac{\dd w_{mn}}{\dd t}=2\epsilon[\alpha(\lambda X\frac{\dd
    W}{\dd t})_i-S_i(\lambda^t X\frac{\dd
    W}{\dd t} \bar{\eta})-\bar{\eta}^t\frac{\dd W}{\dd
    t}\bar{\eta}\lambda_i],$$ 
where
$$L=\bar{\eta}_k\frac{\dd z_k}{\dd t},~\lambda=
\bar{\eta}^t\frac{\dd W}{\dd t},~ Q=\frac{\dd W}{\dd t}X \frac{\dd z}{\dd t}. $$
We write down \eqref{2EE1} as 
\begin{equation}\label{EC1}
G^3=2\epsilon G^1, 
\end{equation}
where
\begin{equation}\label{G1}
G^1_i=(\bar{\eta}^tY)Y_i +(\bar{\eta}^t\frac{\dd W}{\dd t}X
Y)S_i-\alpha(\frac{\dd W}{\dd t}X)Y_i , ~~~Y=\frac{\dd Z}{\dd t}+\frac{\dd
  W}{\dd t}\bar{\eta}, 
\end{equation}
\begin{equation}\label{G3}
G^3=\frac{\dd^2 z}{\dd t^2}+2Q.
\end{equation}


Now we calculate \eqref{2EE2}. We find successively 
$$\Gamma^{pq}_{jk}\frac{\dd z_j}{\dd
  t}\frac{\dd z_k}{\dd t}=2\epsilon\frac{\dd z_p}{\dd t}\frac{\dd z_q}{\dd t},$$
$$\sum_{m\le n}\Gamma^{pq}_{imn}\frac{\dd z_i}{\dd
  t}\frac{\dd w_{mn}}{\dd t}=\epsilon(\lambda_p\frac{\dd z_q}{\dd t}+\lambda_q\frac{\dd z_p}{\dd t }), $$
$$\sum_{m\le n, u\le v}\Gamma^{pq}_{mnuv}\frac{\dd w_{mn}}{\dd
  t}\frac{\dd w_{uv}}{\dd t}=2 \epsilon \lambda_p\lambda_q + 2(\frac{\dd W}{\dd
t} X \frac{\dd W}{\dd t})_{pq} .$$
We write \eqref{2EE2} as 
\begin{equation}\label{EC2}
G^2=-2\epsilon G^4,
\end{equation}
where
\begin{equation}\label{G2G4}
G^2=\frac{\dd^2W}{\dd t^2}+ 2\frac{\dd W}{\dd t}X\frac{\dd W}{\dd
  t}; ~ G^4= Y \otimes Y\end{equation}We rewrite equations \eqref{EC1},
\eqref{EC2} as \eqref{ECG}
\begin{subequations}\label{ECG}
\begin{align}
kG^3 & =2\mu G^1,\label{ECG1}\\
kG^2& =-2\mu G^4,\label{ECG2}
\end{align}
\end{subequations} 
and we have   proved 
\begin{Proposition}\label{PR2}
The equations of geodesics \eqref{2EE}  on the Siegal-Jacobi ball
$\mc{D}^J_n$ 
corresponding to the balanced metric  attached to the K\"ahler two-form
\eqref{aabX} are given in  \eqref{ECG}, where $G^1$, $G^2$ and $G^4$,
$G^3$, are given by \eqref{G1}, \eqref{G2G4}, respectively
\eqref{G3}.

 Equations \eqref{ECG} in the case of $\mc{D}^J_1$ are those given in
 Proposition  \ref{GER}.
\end{Proposition}
Now we discuss the solution of the system of differential equations
\eqref{ECG}.

(a) If we consider $\mu =0$, then in \eqref{ECG2} we get $G^2=0$,
i.e. the equation of geodesics on the Siegel ball $\mc{D}_n$ (see the Appendix)
\begin{equation}\label{geod}
G^2:=\frac{\dd^2W}{\dd t^2}+ 2\frac{\dd W}{\dd t}X\frac{\dd W}{\dd
  t} = 0.\end{equation}
The solution of the equation \eqref{geod} with $W(0)=\db{O}_n$ and
$\dot{W}(0)=B$
is 
\begin{equation}\label{THH}
W(t)=B\frac{\tanh (t\sqrt{\bar{B}B})}{\sqrt{\bar{B}B}}.
\end{equation}
Note that {\it the change of coordinates} \eqref{THH} {\it is a}
$FC$-{\it transform}
in the sense of \eqref{bigch} on the Siegel ball $\mc{D}_n$ (see the
comment after Remark 7 in \cite{csg} for $\mc{D}_1$ and formula (2.25)
for $\mc{D}_n$  in
\cite{sbj}).

If we introduce the  solution \eqref{THH} into \eqref{ECG1}, then the
solution of the equation $G^3=0$ is 
$$z(t)=W(t)B^{-1}\left(\frac{\dd z}{\dd t}\right)_0 +z(0).$$

(b) If $\mu \not=0$, a particular solution $(z,W)$ of the equations of
geodesics on Siegel-ball Jacobi $\mc{D}^J_n$ is given by
$z(t)=\eta_0-W(t)\bar{\eta}_0$, where $W(t)$ has the expression
\eqref{THH} and $\eta_0$ is independent of $t$. This is a particular
case of the solution corresponding to $\eta=\eta_0+t\eta_1$ of the equations of
geodesics $\frac{\dd ^2\eta}{\dd t^2}=0$ on the flat manifold $\C^n$
corresponding to the separation of variables by the $FC$-transform. We
can formulate the
\begin{Proposition}\label{part} The $FC$-transform \eqref{etaZ} is not  a geodesic mapping, but
  it gives geode\-sics $(z(t),W(t))=FC(\eta_0,W(t))$ on the nonsymmetric
  space $\mc{D}^J_n$,
  with $W(t)$ given by \eqref{THH}. 
\end{Proposition} 
For more details, see Remark 8 and appendix A in \cite{csg},
where  the notion of geodesic
  mapping (cf. Definition 5.1 p. 127 in \cite{mikes}) is explained. 
\subsection{Covariant derivative of one-forms}\label{CDDD}
We calculate the connection matrix 
\eqref{CV1} on the Siegel-Jacobi ball using the Christofell's symbols
obtained  in Proposition \ref{CRR}
\begin{equation}\label{mw}
\theta=\left(\begin{array}{cc}\theta^i_j & \theta^i_{pq}\\\theta^{pq}_i&
                                                     \theta^{pq}_{mn}\end{array}\right) .\end{equation}
We have 
\begin{equation}\label{ec1}\theta^i_j:= \Gamma^i_{jk}\dd z_k+\sum_{p\le q}\Gamma^i_{jpq}\dd
w_{pq}.
\end{equation}
We get
\begin{equation}\label{mw1}
\theta^i_j=(1+\epsilon\alpha)\Xi_{ji}-\epsilon[\bar{\eta}_j\mc{A}_i+\delta_{ij}\bar{\eta}^t\mc{A}+S_i(\Xi
\bar{\eta})_j],\end{equation}
where we have introduced the notation 
\begin{equation}\label{XI}\Xi= X\dd W.\end{equation}
We calculate
\begin{equation}\label{ec2}\theta^i_{pq}:=\Gamma^i_{pqj}\dd z_j+\sum_{m\le n}\Gamma^{i}_{pqmn}\dd
w_{mn}, \end{equation}
and we obtain 
\begin{equation}\label{mw2}
\begin{split}
\theta^i_{pq} & =f_{pq}[\delta_{ip}(X\dd z)_q+\delta_{iq}(X\dd
z)_p+\epsilon T^i_{pq}],\\
T^i_{pq} &=
\alpha[\delta_{ip}(X\mc{A})_q+\delta_{iq}(X\mc{A})_p+\bar{\eta}_p\Xi_{qi}+\bar{\eta}_q\Xi_{pi} ]\\
~~ & -S_i\left\{\bar{\eta}_p[X(\mc{A}+\dd W\bar{\eta})]_q+\bar{\eta}_q[X(\mc{A}+\dd W\bar{\eta})]_p\right\}\\
~~& - [(\bar{\eta}_q\delta_{ip}+\bar{\eta}_p\delta_{iq})(\bar{\eta}^t\mc{A})+2\bar{\eta}_p\bar{\eta}_q\mc{A}_i].
\end{split}
\end{equation}
For 
\begin{equation}\label{ec3}\theta^{pq}_i:=\Gamma^{pq}_{ij}\dd z_j+\sum_{m\le n}\Gamma^{pq}_{imn}\dd
w_{mn},\end{equation}
we get
\begin{equation}\label{mw3}
\theta^{pq}_i=
\epsilon(\delta_{iq}\mc{A}_p+\delta_{ip}\mc{A}_q).
\end{equation}
We  calculate 
\begin{equation}\label{ec4}\theta^{pq}_{mn}:=\Gamma^{pq}_{imn}\dd z_i+ \sum_{u\le
  v}\Gamma^{pq}_{mnuv}\dd w_{uv},\end{equation}
as
$$ \Gamma^{pq}_{imn}\dd z_i=\epsilon[\bar{\eta}_m(\delta_{np}\dd
z_q+\delta_{nq}\dd z_p)+\bar{\eta}_n(\delta_{mp}\dd
z_q+\delta_{mq}\dd z_p)]f_{mn},$$
\begin{equation}
\begin{split}
\sum_{u\le
  v}\Gamma^{pq}_{mnuv}\dd w_{uv} & =[\delta_{pm}\Xi_{nq}+
\delta_{pn}\Xi_{mq}+\delta_{qm}\Xi_{np}+\delta_{qn}\Xi_{mp}]f_{mn}\\
& + \epsilon[
(\bar{\eta}_n\delta_{mq}+\bar{\eta}_m\delta_{nq})(\bar{\eta}^t\dd
W)_p+  (\bar{\eta}_n\delta_{mp}+\bar{\eta}_m\delta_{np})(\bar{\eta}^t\dd
W)_q ] f_{mn}
\end{split}
\end{equation}
For \eqref{ec4} we obtain
\begin{equation}\label{mw4}
\begin{split}
\theta^{pq}_{mn} & =[\delta_{pm}\Xi_{nq}+
\delta_{pn}\Xi_{mq}+\delta_{qm}\Xi_{np}+\delta_{qn}\Xi_{mp}]f_{mn}\\
& +
\epsilon[(\bar{\eta}_n\delta_{mq}+\bar{\eta}_m\delta_{nq})\mc{A}_p+(\bar{\eta}_n\delta_{mp}+\bar{\eta}_m\delta_{np})\mc{A}_q ]f_{mn}.
\end{split}
\end{equation}
We calculate the covariant derivative of the one form $\dd z_i$  with \eqref{CV2} as 
\begin{equation}\label{DDZ1}-D (\dd z_i) =w^i_j \dd z_j+ \sum_{p\le q}\theta^i_{pq}
\dd w_{pq}.\end{equation}
We get successively
$$\theta^i_j \dd z_j= (1+\epsilon\alpha)(\Xi^t\dd
z)_i-\epsilon[(\bar{\eta}^t\dd z)\mc{A}_i+(\bar{\eta}^t\mc{A})\dd
z_i+(\bar{\eta}^t\Xi^t\dd z)S_i],$$
\begin{equation}
\begin{split} 
\sum_{p\le q}\theta^i_{pq}\dd w_{pq} & =(\Xi^t\dd
z)_i+\epsilon\{\alpha[\Xi^t(\mc{A}+\dd
W\bar{\eta})]_i+[\bar{\eta}^t\Xi^t(\mc{A}+\dd W\bar{\eta})]S_i\} \\
& -\epsilon[(\bar{\eta}^t\mc{A})(\dd W\bar{\eta})_i+(\bar{\eta}^t\dd W \bar{\eta})\mc{A}_i].
\end{split}
\end{equation}
We get for \eqref{DDZ1} the expression 
\begin{equation}\label{DZI}
-\frac{1}{2}D (\dd z_i)=\{\Xi^t[(1+\epsilon\alpha )\mc{A}-\dd
W\bar{\eta}]\}_i-\epsilon[(\bar{\eta}^t\mc{A})\mc{A}_i+(\bar{\eta}^t\Xi^t\mc{A})S_i]. 
\end{equation}
Now we calculate the covariant derivative of $\dd w_{pq}$:
\begin{equation}\label{DDW2}-D \dd w_{pq}= \theta^{pq}_i\dd z_i +\sum_{m\le n}\theta^{pq}_{mn}\dd w_{mn}.\end{equation}
We obtain
$$\theta^{pq}_i\dd z_i=\epsilon(\mc{A}_p\dd z_q+\mc{A}_q\dd z_p),$$
$$\sum_{m\le n}\theta^{pq}_{mn}\dd w_{mn}=(\dd W\Xi)_{pq}+    (\dd
W\Xi)_{qp}+\epsilon[\mc{A}_p(\dd W \bar{\eta})_q+\mc{A}_q(\dd W
\bar{\eta}) _p] ,   $$
which introduced in \eqref{DDW2} gives
$$-D (\dd w_{pq} )=(\dd W\Xi)_{pq}+    (\dd
W\Xi)_{qp}+2\epsilon\mc{A}_p\mc{A}_q,$$
which can be  writen as
\begin{equation}\label{DW}-D(\dd W)=2\left(\begin{array}{cc}\mc{A}^t \dd
               W\end{array}\right)\left(\begin{array}{cc}\epsilon\un
                                                 & 0\\0 &
                                                          X\end{array}\right)\left(\begin{array}{c}\mc{A}
  \\ \dd W\end{array}\right).
\end{equation}
We have proved
\begin{Proposition}\label{last}
The connections matrix \eqref{mw} on $\mc{D}^J_n$, whose matrix
  elements are calculated with formulae \eqref{ec1}, \eqref{ec2},
  \eqref{ec3} and \eqref{ec4},  are given by
  \eqref{mw1}, \eqref{mw2}, \eqref{mw3}, respectively \eqref{mw4}.

The covariant derivative  $D(\dd z_i)$ \eqref{DDZ1} and  $ D(\dd w_{pq})$
\eqref{DDW2}  on the Siegel-Jacobi ball  $\mc{D}^J_n$ are given by formulae
\eqref{DZI}, respectively \eqref{DW}, which generalize  \eqref{DDZ},
respectively \eqref{DDW} on the Siegl-Jacobi disk $\mc{D}^J_1$.
\end{Proposition}

\section{Appendix}\label{APP}

The real Jacobi group of index $n$ is defined as $G^J_n(\R )=\mr{H}_n(\R)\rtimes\text{Sp}(n,\R
) $, where $\mr{H}_n(\R)$ is the real Heisenberg group
of real dimension $(2n+1)$ \cite{Y02,nou}.

To the Jacobi group  $G^J_n(\R )$ it is associated the homogeneous
manifold  - the Siegel-Jacobi upper half-plane - 
$\mc{X}^J_n\approx\C^n \times
\mc{X}_n$,  where the   Siegel
upper half-plane  $\mc{X}_n= \mr{Sp}(n,\R)/\mr{U}(n)$    is realized as $$\mc{X}_n:=\{V\in M(n,\C)| V=S+\ii  R, S,
R\in M(n,\R), 
R>0,
 S^t=S; R^t=R\} . $$

Siegel has determined the metric on $\mc{X}_n$, 
$\text{Sp}(n,\R)$-invariant to the action \eqref{conf1}
\begin{equation}\label{conf1}
V_1  =  (aV+b)(cV+d)^{-1}
=(Vc^t+d^t)^{-1}(Va^t+b^t); \left(\begin{array}{cc} a & b\\c
                                                      &d\end{array}\right)\in \text{Sp}(n,\R),
\end{equation}
 (see equation (2) in
 \cite{sieg} or Theorem 3 at p. 644 in \cite{hua44}):
\begin{equation}\label{msig}
\dd s^2_{\mc{X}_n}(R)=\tr (R^{-1}\dd V R^{-1} \dd \bar{V}).
\end{equation}
With the Cayley transform  \eqref{big11} 
\begin{equation}\label{big11} 
V  =  \ii (\un-W ) ^{-1} (\un+W ),  
\end{equation}
and the relations 
\begin{subequations}\label{scgimb}
\begin{align}
\dot{V}  & = 2\ii U \dot{W} U, ~~ U=(\un-W)^{-1},\label{scg1}\\
\ddot{V} & = 2\ii U(\ddot{W}+2\dot{W}U\dot{W} )U, \label{scg3}\\
2 R  & =(\un +W)U+(\un+\bar{W})\bar{U} ,\label{scg2}
\end{align}
\end{subequations}
 introduced in \eqref{msig},  it is obtained the metric
on $\mc{D}_n$, $\text{Sp}(n,\R)_{\C}$-invariant to the action
\eqref{x44}: 
\begin{equation}\label{mtrball}
\dd s^2_{\mc{D}_n}(W)= 4\tr(M\dd
W\bar{M}\dd \bar{W}),~ W\in\mc{D}_n, ~M=(\un-W\bar{W})^{-1}.
\end{equation}
\eqref{scg3}  is equation (3.2) in \cite{Y10} and corresponds to
$\omega_{\mc{D}_n}$ in \eqref{2omega}.

Equations of geodesics on the Siegel-upper half plane $\mc{X}_n$ (see
e.g. equation (39) in \cite{sieg} or Theorem 11 in at
p. 478 in \cite{hua44}) are
\begin{equation}\label{geoS}\ddot{V}=-\ii
  \dot{V}R^{-1}\dot{V}.\end{equation}

The  relation \eqref{scg2}  can be written as
\begin{equation}\label{Yr}
R= (\un-\bar{W})^{-1}(\un-\bar{W}W)(\un-W)^{-1}= (\un-W)^{-1}(\un-W\bar{W})(\un-\bar{W})^{-1}.
\end{equation}

Equation \eqref{geoS}  on $\mc{X}_n$ becomes on $\mc{D}_n$
$$\ddot{W}+2\dot{W}\Sigma\dot{W}=0$$
where
\begin{equation}
\begin{split}
\Sigma  & =\tau U,\\
\tau &= \un -[v(\un-W)]^{-1}\\
 ~~ &  =2\bar{U}(\un-\bar{W}W),
\end{split}
\end{equation}
and we have $$\Sigma = \bar{W}(\un-W\bar{W})^{-1},$$
i.e. $\Sigma$ is our  previous $X$ defined in \eqref{LIK} and
  equation \eqref{geod} is retrieved. Note that the equation of
  geodesics on the Siegel ball $\mc{D}_n$
has the same expression  as on the noncompact Grassmannian $\text{SU}(n,n)/\text{S}(\text{U}(n)\times\text{U}(n))$,
see equation (6.13)  in \cite{gras}.

{\bf Acknowledgement} The author is grateful to Mihai Visinescu
 for clarifying  discussions on connections. This research  was conducted in  the  framework of the 
ANCS project  program   PN 16 42 01 01/2016 and   UEFISCDI-Romania 
 program PN-II-PCE-55/05.10.2011.

\end{document}